\documentstyle{amsppt}
\magnification1095
\voffset-7.5mm
\hoffset-2.5mm
\pagewidth{38.33pc} \pageheight{52.5pc}
\baselineskip=12pt plus0.2pt minus0.2pt
\abovedisplayskip=7pt plus2pt minus2pt
\belowdisplayskip=7pt plus2pt minus2pt
\abovedisplayshortskip=5pt plus1pt minus3pt
\belowdisplayshortskip=5pt plus1pt minus3pt
\hfuzz=2.5pt\rightskip=0pt plus1pt
\binoppenalty=10000\relpenalty=10000\relax
\TagsOnRight
\loadbold
\nologo
\addto\tenpoint{\normalbaselineskip=1.2\normalbaselineskip\normalbaselines}
\addto\eightpoint{\normalbaselineskip=1.2\normalbaselineskip\normalbaselines}
\let\wt\tilde
\define\dn{\vphantom{{}_0}}
\topmatter
\title
A third-order Ap\'ery-like recursion for $\zeta(5)$
\endtitle
\author
Wadim Zudilin \rm(Moscow)
\endauthor
\date
\hbox to100mm{\vbox{\hsize=100mm%
\centerline{E-print \tt math.NT/0206178}
\smallskip
\centerline{27 May 2002}
\centerline{{\sl Last revision\/}: 06 July 2002}
}}
\enddate
\address
\hbox to70mm{\vbox{\hsize=70mm%
\leftline{Moscow Lomonosov State University}
\leftline{Department of Mechanics and Mathematics}
\leftline{Vorobiovy Gory, GSP-2, 119992 Moscow, RUSSIA}
\leftline{{\it URL\/}: \tt http://wain.mi.ras.ru/index.html}
}}
\endaddress
\email
{\tt wadim\@ips.ras.ru}
\endemail
\endtopmatter
\rightheadtext{A third-order Ap\'ery-like recursion for $\zeta(5)$}
\leftheadtext{W.~Zudilin}
\document

In 1978, R.~Ap\'ery~\cite{1},~\cite{2} has given sequences
of rational approximations to~$\zeta(2)$ and $\zeta(3)$ yielding
the irrationality of each of these numbers. One of the key
ingredient of Ap\'ery's proof are second-order difference equations
with polynomial coefficients satisfied by numerators
and denominators of the above approximations.
Recently, V.\,N.~Sorokin~\cite{3} and this author~\cite{4},~\cite{5}
have got independently a similar second-order difference equation
for~$\zeta(4)$. This recursion does not give diophantine approximations
to~$\zeta(4)=\pi^4/90$ proving its irrationality,
however it presents an algorithm
for fast calculation of this constant. The aim of this note
is to highlight a possible generalization of the above results
for the number~$\zeta(5)$, the irrationality of which is not
yet proved since nowadays.

\subhead
1. Statement of the main result
\endsubhead
Consider the difference equation
$$
(n+1)^6a_0(n)q_{n+1}+a_1(n)q_n-4(2n-1)a_2(n)q_{n-1}
-4(n-1)^4(2n-1)(2n-3)a_0(n+1)q_{n-2}=0,
\tag1
$$
where
$$
\align
a_0(n)
&=41218n^3-48459n^2+20010n-2871,
\\
a_1(n)
&=2(48802112n^9+89030880n^8+36002654n^7-24317344n^6-19538418n^5
\\ &\qquad
+1311365n^4+3790503n^3+460056n^2-271701n-60291),
\\
a_2(n)
&=3874492n^8-2617900n^7-3144314n^6+2947148n^5+647130n^4
\\ &\qquad
-1182926n^3+115771n^2+170716n-44541,
\endalign
$$
and define its three linearly independent solutions
$\{q_n\}$, $\{p_n\}$, and $\{\wt p_n\}$ by the initial data
$$
\gathered
q_0=-1, \quad q_1=42, \quad q_2=-17934,
\qquad
p_0=0, \quad p_1=\frac{87}2, \quad p_2=-\frac{1190161}{64},
\\
\wt p_0=0, \quad \wt p_1=\frac{101}2,
\quad \wt p_2=-\frac{344923}{16}
\endgathered
$$
(here and below a symbol $\{x_n\}$ means the sequence
$\{x_n\}_{n=0}^\infty=\{x_0,x_1,x_2,\dots\}$).

\proclaim{Theorem 1}
The sequences
$$
\ell_n=q_n\zeta(5)-p_n, \quad \wt\ell_n=q_n\zeta(3)-\wt p_n,
\qquad n=0,1,2,\dots,
\tag2
$$
that also satisfy the difference equation~\thetag{1},
are of constant sign:
$$
\ell_n>0, \quad \wt\ell_n<0,
\qquad n=1,2,\dots,
\tag3
$$
and the following limit relations hold:
$$
\gathered
\lim_{n\to\infty}\frac{\log|\ell_n|}n
=\lim_{n\to\infty}\frac{\log|\wt\ell_n|}n
=\log|\mu_2|=-1.08607936\dots,
\\
\lim_{n\to\infty}\frac{\log|q_n|}n
=\lim_{n\to\infty}\frac{\log|p_n|}n
=\lim_{n\to\infty}\frac{\log|\wt p_n|}n
=\log|\mu_3|,
\endgathered
\tag4
$$
where
$$
\mu_1=-0.02001512\dots, \qquad
\mu_2=0.33753726\dots, \qquad
\mu_3=-2368.31752213\dots
\tag5
$$
are the roots of the characteristic polynomial
$\mu^3+2368\mu^2-752\mu-16$ of the recursion~\thetag{1}.
\endproclaim

Note that the sequences
$\{q_n\}$, $\{p_n\}$, and $\{\wt p_n\}$
are of alternating sign:
$$
(-1)^{n-1}q_n>0, \quad (-1)^{n-1}p_n>0, \quad (-1)^{n-1}\wt p_n>0,
\qquad n=1,2,\dots\,.
$$

Theorem~1 gives an algorithm for fast calculation of
the number~$\zeta(5)$. Namely, the sequence of rational
numbers $p_n/q_n$ converges to~$\zeta(5)$ with speed
$|\mu_2/\mu_3|<1.42521964\cdot10^{-4}$
(see the table).

\medskip
\line{\hss\vbox{\offinterlineskip
\halign to85.5mm{\strut\tabskip=1mm minus 1mm
\strut\vrule#&\hbox to7mm{\hss$#$\hss}&%
\vrule#&\hbox to49mm{\hss$#$\hss}&%
\vrule#&\hbox to29mm{\hss$#$\hss}&%
\vrule#\tabskip=0pt\cr\noalign{\hrule}
height9.2pt&n&&p_n/q_n&&|\zeta(5)-p_n/q_n|&\cr
\noalign{\hrule\vskip1pt\hrule}
height9.2pt&0&&0&&1.036927755\dots&\cr
\noalign{\hrule}
height10.5pt&1&&\frac{29}{28\dn}&&0.001213469\dots&\cr
\noalign{\hrule}
height10.5pt&2&&\frac{24289}{23424\dn}&&0.000000182\dots&\cr
\noalign{\hrule}
height10.5pt&3&&\frac{7682021239}{7408444032\dn}&&<2.80\cdot10^{-11}\;\,&\cr
\noalign{\hrule}
height10.5pt&4&&\frac{24943788950905}{24055474286592\dn}&%
&<4.13\cdot10^{-15}\;\,&\cr
\noalign{\hrule}
height10.5pt&5&&\frac{81875586674776013003}{78959779279372800000\dn}&%
&<6.02\cdot10^{-19}\;\,&\cr
\noalign{\hrule}
height10.5pt&6&&\frac{282653756112686336975107}{272587704119854963200000\dn}&%
&<8.71\cdot10^{-23}\;\,&\cr
\noalign{\hrule}
height10.5pt&7&%
&\frac{215903781003833520407770175189}{208214873150908926517286400000\dn}&%
&<1.26\cdot10^{-26}\;\,&\cr
\noalign{\hrule}
height10.5pt&10&& &&<3.71\cdot10^{-38}\;\,&\cr
\noalign{\hrule}
height10.5pt&20&& &&<1.32\cdot10^{-76}\;\,&\cr
\noalign{\hrule}
height10.5pt&50&& &&<5.52\cdot10^{-192}&\cr
\noalign{\hrule}
}}\hss}

\smallskip
Further construction of the above solutions of the
difference equation~\thetag{1} leads to the inclusions
$$
4D_n^2q_n\in\Bbb Z, \quad
4D_n^7p_n\in\Bbb Z, \quad
4D_n^5\wt p_n\in\Bbb Z, \qquad
n=1,2,\dots
$$
(see relations~\thetag{14} and~\thetag{15}),
where $D_n$~denotes the least common multiple of numbers
$1,2,\dots,n$, while straightforward calculations on the basis of
the recursion~\thetag{1} show that
$$
q_n\in\Bbb Z, \quad
2D_n^5p_n\in\Bbb Z, \quad
2D_n^3\wt p_n\in\Bbb Z, \qquad
n=1,2,\dots\,.
\tag6
$$
The inclusions~\thetag{6} do not allow oneself to prove
the irrationality of the number~$\zeta(5)$, hence we do not accent
our attention on arithmetic properties of the sequences
$\{q_n\}$, $\{p_n\}$, $\{\wt p_n\}$
and stress only on an applied meaning of Theorem~1.

\subhead
2. Auxiliary recursions
\endsubhead
The very-well-poised hypergeometric series
$$
\aligned
r_n
&=\phantom-n!^4\sum_{k=1}^\infty\Bigl(k+\frac n2\Bigr)
\frac{\prod_{j=1}^n(k-j)\cdot\prod_{j=1}^n(k+n+j)}
{\prod_{j=0}^n(k+j)^6},
\\
\wt r_n
&=-n!^4\sum_{k=1}^\infty\Bigl(k+\frac n2\Bigr)
\frac{\prod_{j=0}^n(k-j)\cdot\prod_{j=0}^n(k+n+j)}
{\prod_{j=0}^n(k+j)^6},
\endaligned
\qquad n=0,1,2,\dots,
\tag7
$$
are $\Bbb Q$-linear forms in $1,\zeta(3),\zeta(5)$:
$$
r_n=u_n\zeta(5)+w_n\zeta(3)-v_n,
\quad
\wt r_n=\wt u_n\zeta(5)+\wt w_n\zeta(3)-\wt v_n,
\qquad n=0,1,2,\dots\,.
$$
An easy verification shows that
$$
\alignedat3
r_0&=\zeta(5), \qquad&
r_1&=9\zeta(5)+33\zeta(3)-49, \qquad&
r_2&=469\zeta(5)+\frac{6125}4\zeta(3)-\frac{74463}{32},
\\
\wt r_0&=\zeta(3), \qquad&
\wt r_1&=2\zeta(5)+12\zeta(3)-\frac{33}2, \qquad&
\wt r_2&=552\zeta(5)+1764\zeta(3)-\frac{43085}{16}.
\endalignedat
$$
Applying the algorithm of creative telescoping~\cite{6, Chapter~6}
to the series~\thetag{7} in a manner
of the works~\cite{4}, \cite{5}, \cite{7},
we arrive at the difference equations
$$
n(n+1)^5b_0(n-1)u_{n+1}-2nb_1(n)u_n-b_2(n)u_{n-1}
+2(n-1)^5(2n-1)b_0(n)u_{n-2}=0,
\tag8
$$
where
$$
b_0(n)=-a_0(-n),
\qquad
b_1(n)=a_2(-n),
\qquad
b_2(n)=-a_1(-n),
$$
and
$$
n^3(n+1)^3\wt b_0(n-1)u_{n+1}-2n\wt b_1(n)u_n-\wt b_2(n)u_{n-1}
+2n(n-1)^4(2n-3)\wt b_0(n)\wt u_{n-2}=0,
\tag9
$$
where
$$
\align
\wt b_0(n)
&=41218n^7+35648n^6-932n^5-13190n^4-5128n^3+811n^2+957n+174,
\\
\wt b_1(n)
&=3874492n^{12}-14084302n^{11}+12425954n^{10}+8641603n^9
-15230839n^8-1369195n^7
\\ &\qquad
+8618417n^6-623249n^5-2785973n^4
+308165n^3+495325n^2-40670n-37632,
\\
\wt b_2(n)
&=2(48802112n^{13}-201803328n^{12}+267014032n^{11}-69927236n^{10}
\\ &\qquad
-95912858n^9+37524471n^8+30257812n^7-9523224n^6-8524312n^5
\\ &\qquad
+2138687n^4+1507490n^3-398634n^2-111012n+33408).
\endalign
$$
Note that each of the recursions~\thetag{8} and~\thetag{9}
has the same characteristic polynomial
$\lambda^3-188\lambda^2-2368\lambda+4$; its roots
$\lambda_1,\lambda_2,\lambda_3$ ordered in increasing order
of their moduli are related to the roots~\thetag{5}
as follows:
$$
\mu_1=\lambda_1\lambda_2, \qquad
\mu_2=\lambda_1\lambda_3, \qquad
\mu_3=\lambda_2\lambda_3.
$$

\proclaim{Theorem 2}
The sequences of the linear forms $\{r_n\}$
and their coefficients $\{u_n\}$, $\{w_n\}$, $\{v_n\}$
satisfy the difference equation~\thetag{8}.
In addition, the inequalities
$$
r_n>0, \quad u_n>0, \quad w_n>0, \quad v_n>0,
\qquad n=1,2,\dots,
\tag10
$$
and the limit relations
$$
\lim_{n\to\infty}\frac{\log|r_n|}n=\log\lambda_1,
\qquad
\lim_{n\to\infty}\frac{\log|u_n|}n
=\lim_{n\to\infty}\frac{\log|w_n|}n
=\lim_{n\to\infty}\frac{\log|v_n|}n
=\log\lambda_3
\tag11
$$
hold.
\endproclaim

\proclaim{Theorem 3}
The sequences of the linear forms $\{\wt r_n\}$
and their coefficients $\{\wt u_n\}$, $\{\wt w_n\}$, $\{\wt v_n\}$
satisfy the difference equation~\thetag{9}.
In addition, the inequalities
$$
\wt r_n<0, \quad \wt u_n>0, \quad \wt w_n>0, \quad \wt v_n>0,
\qquad n=1,2,\dots,
\tag12
$$
and the limit relations
$$
\lim_{n\to\infty}\frac{\log|\wt r_n|}n=\log\lambda_1,
\qquad
\lim_{n\to\infty}\frac{\log|\wt u_n|}n
=\lim_{n\to\infty}\frac{\log|\wt w_n|}n
=\lim_{n\to\infty}\frac{\log|\wt v_n|}n
=\log\lambda_3
\tag13
$$
hold.
\endproclaim

We mention that Vasilyev's result~\cite{8} and our theorem~\cite{4}
(see also \cite{9})
on the coincidence of very-well-poised hypergeometric series
and a suitable generalization of Beukers' integral~\cite{10}
lead one to the inclusions
$$
\gathered
2u_n\in\Bbb Z, \quad 2D_n^2w_n\in\Bbb Z,
\quad 2D_n^5v_n\in\Bbb Z,
\\
2\wt u_n\in\Bbb Z, \quad 2D_n^2\wt w_n\in\Bbb Z,
\quad 2D_n^5\wt v_n\in\Bbb Z,
\endgathered
\qquad n=1,2,\dots\,.
\tag14
$$

Finally, the required sequences~\thetag{2} are determined by the formulae
$$
\aligned
\ell_n
&=\wt w_nr_n-w_n\wt r_n
=(u_n\wt w_n-\wt u_nw_n)\zeta(5)
-(\wt w_nv_n-w_n\wt v_n),
\\
\wt\ell_n
&=u_n\wt r_n-\wt u_nr_n
=(u_n\wt w_n-\wt u_nw_n)\zeta(3)
-(u_n\wt v_n-\wt u_nv_n),
\endaligned
\qquad n=0,1,2,\dots,
$$
hence
$$
q_n=u_n\wt w_n-\wt u_nw_n, \quad
p_n=\wt w_nv_n-w_n\wt v_n, \quad
\wt p_n=u_n\wt v_n-\wt u_nv_n,
\qquad n=0,1,2,\dots\,.
\tag15
$$
For proving the inequalities~\thetag{3} and the limit
relations~\thetag{4}, it remains to use the estimates
\thetag{10},~\thetag{12},
relations~\thetag{11},~\thetag{13}, and
Poincar\'e's theorem.

\subhead
3. Conclusion remarks
\endsubhead
Recently, another algorithm
for fast calculation of~$\zeta(5)$ and~$\zeta(3)$
(also not producing good enough diophantine
approximations to these constants),
based on the infinite matrix product
$$
\prod_{n=1}^\infty\pmatrix
-\dfrac n{2(2n+1)} & \dfrac1{2n(2n+1)} & \dfrac1{n^4} \\
0 & -\dfrac n{2(2n+1)} & \dfrac5{4n^2} \\ \vspace{3pt}
0 & 0 & 1 \endpmatrix
=\pmatrix 0 & 0 & \zeta(5) \\ 0 & 0 & \zeta(3) \\ 0 & 0 & 1 \endpmatrix,
$$
has been suggested by R.\,W.~Gosper~\cite{11}. However the letters~\cite{11}
do not contain any description of analytic and arithmetic properties
of the rational approximations so constructed.
Mention also an algorithm due to E.\,A.~Karatsuba~\cite{12}
for calculation of the Riemann zeta function
at positive integers.

The above scheme of Section~2 also allows one to construct
a third-order recursion for simultaneous rational
approximations to~$\zeta(2)$ and~$\zeta(3)$.
Namely, consider the difference equation
$$
(n+1)^4\wt a_0(n)q_{n+1}-\wt a_1(n)q_n+4(2n-1)\wt a_2(n)q_{n-1}
-4(n-1)^2(2n-1)(2n-3)\wt a_0(n+1)q_{n-2}=0,
\tag16
$$
where
$$
\align
\wt a_0(n)
&=946n^2-731n+153,
\\
\wt a_1(n)
&=2(104060n^6+127710n^5+12788n^4-34525n^3-8482n^2+3298n+1071),
\\
\wt a_2(n)
&=3784n^5-1032n^4-1925n^3+853n^2+328n-184,
\endalign
$$
and define its three linearly independent solutions
$\{q_n'\}$, $\{p_n'\}$, and $\{\wt p_n'\}$ by the initial data
$$
q_0'=1, \quad q_1'=14, \quad q_2'=978,
\qquad
p_0'=0, \quad p_1'=17, \quad p_2'=\frac{9405}8,
\qquad
\wt p_0'=0, \quad \wt p_1'=23, \quad \wt p_2'=\frac{6435}4.
$$

\proclaim{Theorem 4}
The sequences $\{q_n'\}$, $\{p_n'\}$, and $\{\wt p_n'\}$
\rom(of positive sign\rom) as well as the sequences
$$
\ell_n'=q_n'\zeta(3)-p_n', \quad \wt\ell_n'=q_n'\zeta(2)-\wt p_n',
\qquad n=0,1,2,\dots,
$$
satisfy the limit relations
$$
\gather
\lim_{n\to\infty}\frac{\log|\ell_n'|}n
=\lim_{n\to\infty}\frac{\log|\wt\ell_n'|}n
=\log|\mu_2|=-1.31018925\dots,
\\
\lim_{n\to\infty}\frac{\log|q_n'|}n
=\lim_{n\to\infty}\frac{\log|p_n'|}n
=\lim_{n\to\infty}\frac{\log|\wt p_n'|}n
=\log|\mu_3|,
\endgather
$$
where
$$
\mu_{1,2}=0.07260980\ldots\pm i0.25981363\dots, \qquad
\mu_3=219.85478039\dots
$$
are the roots of the characteristic polynomial
$\mu^3-220\mu^2+32\mu-16$ of the recursion~\thetag{16}.
\endproclaim

This time, auxiliary recursions are satisfied
by the hypergeometric series
$$
r_n'=-n!^2\sum_{k=1}^\infty
\frac{\prod_{j=1}^n(k-j)}{\prod_{j=0}^n(k+j)^3},
\quad
\wt r_n'=n!^2\sum_{k=1}^\infty
\frac{\prod_{j=0}^n(k-j)}{\prod_{j=0}^n(k+j)^3},
\qquad n=0,1,2,\dots,
$$
that are $\Bbb Q$-linear forms in $1$, $\zeta(2)$, and $\zeta(3)$.

Similarly to the case of the difference equation~\thetag{1},
the explicit calculations lead to the inclusions
$$
q_n'\in\Bbb Z, \quad D_n^3p_n'\in\Bbb Z, \quad D_n^2\wt p_n'\in\Bbb Z,
\qquad n=1,2,\dots,
$$
that are quite better than one can expect.



\Refs

\ref\no1
\by R.~Ap\'ery
\paper Irrationalit\'e de~$\zeta(2)$ et~$\zeta(3)$
\jour Ast\'erisque
\vol61
\yr1979
\pages11--13
\endref

\ref\no2
\by A.~van der Poorten
\paper A proof that Euler missed...
Ap\'ery's proof of the irrationality of~$\zeta(3)$
\paperinfo An informal report
\jour Math. Intelligencer
\vol1
\issue4
\yr1978/79
\pages195--203
\endref

\ref\no3
\by V.\,N.~Sorokin
\book One algorithm for fast calculation of~$\zeta(4)$
\bookinfo Preprint (April 2002)
\publaddr Moscow
\publ Russian Academy of Sciences,
M.\,V.~Kel\-dysh Institute for Applied Mathematics 
\yr2002
\endref

\ref\no4
\by W.~Zudilin
\paper Well-poised hypergeometric service
for diophantine problems of zeta values
\paperinfo Actes des 12\`emes rencontres arithm\'e\-tiques de Caen
(June 29--30, 2001)
\jour J. Th\'eorie Nombres Bordeaux
\yr2003
\toappear
\endref

\ref\no5
\by W.~Zudilin
\paper Ap\'ery-like difference equation for Catalan's constant
\inbook E-print {\tt math.NT/\allowlinebreak 0201024} (January 2002)
\endref

\ref\no6
\by M.~Petkov\v sek, H.\,S.~Wilf, and D.~Zeilberger
\book $A=B$
\publaddr Wellesley (M.A.)
\publ A.\,K.~Peters, Ltd.
\yr1997
\endref

\ref\no7
\by W.~Zudilin
\paper An elementary proof of Ap\'ery's theorem
\inbook E-print {\tt math.NT/0202159} (February 2002)
\endref

\ref\no8
\by D.\,V.~Vasilyev
\book On small linear forms for the values of the Riemann
zeta-function at odd points
\bookinfo Preprint no.~1\,(558)
\publ Nat. Acad. Sci. Belarus, Institute Math.
\publaddr Minsk
\yr2001
\endref

\ref\no9
\by W.~Zudilin
\paper Very-well-poised hypergeometric series and multiple integrals
\jour Uspekhi Mat. Nauk [Russian Math. Surveys]
\vol57
\yr2002
\issue4
\moreref
\inbook E-print {\tt math.CA/0206177} (March 2002)
\endref

\ref\no10
\by F.~Beukers
\paper A note on the irrationality of~$\zeta(2)$ and~$\zeta(3)$
\jour Bull. London Math. Soc.
\vol11
\yr1979
\issue3
\pages268--272
\endref

\ref\no11
\by R.\,W.~Gosper
\paper Some infinite products involving $\zeta(2m+1)$
\paperinfo Letters of 24.10.2000 and 25.10.2000
\inbook Favorite Mathematical Constants
\ed S.~Finch
\publaddr {\tt http://pauillac.inria.fr/algo/%
bsolve/constant/apery/infprd.html}
\yr2000
\endref

\ref\no12
\by E.\,A.~Karatsuba
\paper Fast computation of the Riemann zeta function~$\zeta(s)$
for integer values of~$s$
\jour Problemy Peredachi Informatsii [Problems Inform. Transmission]
\vol31
\issue4
\yr1995
\pages69--80
\endref

\endRefs
\enddocument
\end